\documentclass[12pt,reqno]{amsart}
\usepackage{amssymb}
\input amssym.def
\usepackage{amsmath,amsfonts,hyperref,xcolor,mathtools,booktabs}
\usepackage{amscd}
\usepackage[mathscr]{eucal}
\usepackage{enumitem}
\usepackage{orcidlink}
\allowdisplaybreaks
\numberwithin{equation}{section}

\hypersetup{colorlinks=true,citecolor={purple},linkcolor={teal},urlcolor={violet}}

\theoremstyle{plain}
\newtheorem{theorem}{Theorem}[section]

\newtheorem{lemma}[theorem]{Lemma}

\theoremstyle{definition}

\makeatletter
\@namedef{subjclassname@2020}{%
	\textup{2020} Mathematics Subject Classification}
\makeatother

\begin{document}
\title[Congruences for an analogue of Lin's partition function]{Congruences for an analogue of Lin's partition function}
\author[Russelle Guadalupe]{Russelle Guadalupe\orcidlink{0009-0001-8974-4502}}
\address{Institute of Mathematics, University of the Philippines Diliman\\
Quezon City 1101, Philippines}
\email{rguadalupe@math.upd.edu.ph}

\renewcommand{\thefootnote}{}



\renewcommand{\thefootnote}{\arabic{footnote}}

\setcounter{footnote}{0}

\begin{abstract}
We study certain arithmetic properties of an analogue $B(n)$ of Lin's restricted partition function that counts the number of partition triples $\pi=(\pi_1,\pi_2,\pi_3)$ of $n$ such that $\pi_1$ and $\pi_2$ comprise distinct odd parts and $\pi_3$ consists of parts divisible by $4$. With the help of elementary $q$-series techniques and modular functions, we establish Ramanujan-type congruences modulo $2,3,5,7$, and $9$ for certain sums involving $B(n)$.
\end{abstract}

\maketitle

\section{Introduction}\label{sec1}

We denote $f_m := \prod_{n\geq 1}(1-q^{mn})$ for $m\in\mathbb{N}$ and $q\in\mathbb{C}$ with $|q| < 1$ throughout this paper. Recall that a (unrestricted) partition of a positive integer $n$ is a finite nonincreasing sequence of positive integers whose sum is $n$. These integers are referred to as the parts of $n$. The generating function for the number $p(n)$ of partitions of $n$ is given by 
\begin{align*}
\sum_{n=0}^\infty p(n)q^n = \dfrac{1}{f_1}.
\end{align*}
Ramanujan \cite[pp. 210--213]{ramp} found the remarkable identities 
\begin{align}
\sum_{n=0}^\infty p(5n+4)q^n &= 5\dfrac{f_5^5}{f_1^6},\label{eq11}\\
\sum_{n=0}^\infty p(7n+5)q^n &= 7\dfrac{f_7^3}{f_1^4}+49q\dfrac{f_7^7}{f_1^8},\nonumber
\end{align}
which yield the congruences  
\begin{align*}
p(5n+4) &\equiv 0\pmod{5},\\
p(7n+5) &\equiv 0\pmod{7}
\end{align*}
for all $n\geq 0$. In 2010, Chan \cite{chan} introduced a cubic partition of $n$, which is a partition of $n$ whose even parts may appear in one of two colors. The generating function for the number $a(n)$ of cubic partitions of $n$ is given by 
\begin{align*}
\sum_{n=0}^\infty a(n)q^n = \dfrac{1}{f_1f_2}.
\end{align*}
Chan \cite{chan} used a $3$-dissection formula involving the Ramanujan's cubic continued fraction to deduce that
\begin{align}\label{eq12}
\sum_{n=0}^\infty a(3n+2)q^n = 3\dfrac{f_3^3f_6^3}{f_1^4f_2^4},
\end{align}
which is an analogue of (\ref{eq11}) and follows that
\begin{align*}
a(3n+2)\equiv 0\pmod{3}
\end{align*}
for all $n\geq 0$. Kim \cite{kim} defined an overcubic partition function $\overline{a}(n)$ of $n$ that counts the number of cubic partitions of $n$ in which the first occurrence of each part may be overlined, whose generating function is 
\begin{align*}
\sum_{n= 0}^\infty \overline{a}(n)q^n = \dfrac{f_4}{f_1^2f_2}.
\end{align*}
Kim \cite{kim} established the following identity
\begin{align}\label{eq13}
\sum_{n= 0}^\infty \overline{a}(3n+2)q^n = 6\dfrac{f_3^6f_4^3}{f_1^8f_2^3}
\end{align}
similar to (\ref{eq11}) and (\ref{eq12}) by employing modular forms, which implies that
\begin{align*}
\overline{a}(3n+2)\equiv 0\pmod{6}
\end{align*}
for all $n\geq 0$. Subsequently, Hirschhorn \cite{hirsc} gave an elementary proof of (\ref{eq13}) and derived the generating functions for $\overline{a}(3n)$ and $\overline{a}(3n+1)$.

In 2013, Lin \cite{lin} studied arithmetic properties of the restricted partition function $b(n)$ that counts the number of partition triples of $\pi=(\pi_1,\pi_2,\pi_3)$ of $n$ such that $\pi_1$ consists of distinct odd parts, and $\pi_2$ and $\pi_3$ consist of parts divisible by $4$. The generating function for $b(n)$ is given by 
\begin{align*}
\sum_{n= 0}^\infty b(n)q^n = \dfrac{f_2^2}{f_1f_4^3}.
\end{align*}
Lin used modular forms to show that
\begin{align}
\sum_{n\geq 0} b(3n+2)q^n &= 3q\dfrac{f_2^6f_{12}^6}{f_1^3f_4^{11}},\label{eq14}\\
\sum_{n\geq 0} b(3n+1)q^n &= \alpha(q^4)\dfrac{f_2^6f_{12}^3}{f_1^3f_4^{10}},\label{eq15}
\end{align}
where $\alpha(q)$ is the cubic theta function of Borwein, Borwein, and Garvan \cite{bbg} defined by
\begin{align*}
\alpha(q) := \sum_{m, n=-\infty}^{\infty} q^{m^2+mn+n^2}.
\end{align*}
We remark that (\ref{eq14}) is an analogue of (\ref{eq12}) and (\ref{eq13}), and yields the congruence 
\begin{align*}
b(3n+2)\equiv 0\pmod{3}
\end{align*}
for all $n\geq 0$. Recently, the author \cite{guad1} provided elementary proofs of (\ref{eq14}) and (\ref{eq15}), and proved certain families of internal congruences modulo $3$ for $b(n)$. 

The objective of this paper is to explore arithmetic properties of a new restricted partition function $B(n)$ analogous to $b(n)$, which counts the number of partition triples $\pi=(\pi_1,\pi_2,\pi_3)$ of $n$ such that $\pi_1$ and $\pi_2$ comprise distinct odd parts, and $\pi_3$ consists of parts divisible by $4$. The generating function for $B(n)$ is then given by 
\begin{align*}
\sum_{n=0}^\infty B(n)q^n = \dfrac{f_2^4}{f_1^2f_4^3}.
\end{align*}
Specifically, we apply elementary $q$-series techniques and modular functions to derive congruences for $B(n)$. Our first two results show congruences modulo $2,3$, and $5$ for $B(n)$.

\begin{theorem}\label{thm11}
For all $n\geq 0$, we have
\begin{align}
	B(2n+1)&\equiv 0\pmod{2},\label{eq16}\\
	B(5n+4)&\equiv 0\pmod{5}.\label{eq17}
\end{align}
\end{theorem}

\begin{theorem}\label{thm12}
For all $n\geq 0$, we have $B(27n+16)\equiv 0\pmod{3}$.
\end{theorem}

The next results reveal congruences modulo $3,7$, and $9$ for certain finite sums involving $B(n)$. Congruences for finite sums involving other restricted partition functions were demonstrated by several authors \cite{amdm,balm,baru,bis1,bis2,guad2,hyz}.

\begin{theorem}\label{thm13}
For all $n\geq 0$ and $j\in \{1,2\}$, we have 
\begin{align*}
	\sum_{k=-\infty}^\infty (-1)^kB\left(9n+3j+2-6k(3k+1)\right)\equiv 0\pmod{3}.
\end{align*}
\end{theorem}

\begin{theorem}\label{thm14}
For all primes $p\equiv 3\pmod{4}$ with $p\geq 5$, $n\geq 0$, and $r\in\{1,\ldots,p-1\}$, we have 
\begin{align*}	
	\sum_{k=-\infty}^\infty (-1)^kB\left(9p^2n+9pr+\dfrac{9(p^2-1)}{4}+2-6k(3k+1)\right)\equiv 0\pmod{3}.
\end{align*}
\end{theorem}

\begin{theorem}\label{thm15}
For all primes $p\equiv 7, 11\pmod{12}$, $n\geq 0$, and $r\in\{1,\ldots,p-1\}$, we have 	
\begin{align*}	
	\sum_{k=-\infty}^\infty (-1)^kB\left(3p^2n+3pr+\dfrac{5(p^2-1)}{4}+1-6k(3k+1)\right)\equiv 0\pmod{9}.
\end{align*}
\end{theorem}

\begin{theorem}\label{thm16}
For all $n\geq 0$, we have 
\begin{align*}
	\sum_{k=-\infty}^\infty (-1)^k(3k+1)B\left(81n+70-54k(3k+2)\right)\equiv 0\pmod{9}.
\end{align*}
\end{theorem}

\begin{theorem}\label{thm17}
For all $n\geq 0$ and $j\in\{3,4,6\}$, we have 
\begin{align}
	\sum_{k=-\infty}^\infty (6k+1)B\left(49n+7j+2-7k(3k+1)\right)\equiv 0\pmod{7},\label{eq18}\\
	\sum_{k=-\infty}^\infty (6k+1)B\left(343n+49j+16-7k(3k+1)\right)\equiv 0\pmod{7}.\label{eq19}
\end{align}
\end{theorem}

We organize the remainder of the paper as follows. We employ classical $q$-series manipulations, $q$-series identities, and dissection formulas to establish Theorems \ref{thm11} in Section \ref{sec2}, Theorems \ref{thm13} and \ref{thm14} in Section \ref{sec3}, and Theorems \ref{thm12}, \ref{thm15}, and \ref{thm16} in Section \ref{sec4}. In particular, we exhibit the exact generating function for $B(3n+2)$ in Section \ref{sec3} and the generating functions modulo $9$ for $B(3n+1)$ and $B(9n+7)$ in Section \ref{sec4} to deduce Theorems \ref{thm12}--\ref{thm16}. In Section \ref{sec5}, we rely on modular functions, particularly the implementation of Radu's Ramanujan--Kolberg algorithm \cite{radu} due to Smoot \cite{smoot}, to prove Theorem \ref{thm17} by finding the generating function modulo $7$ for $B(7n+2)$. 

In the proofs of our main results, we have extensively used without further notice that
\begin{align*}
f_m^{p^k}\equiv f_{mp}^{p^{k-1}}\pmod{p^k}
\end{align*}
for all primes $p$ and integers $m, k\geq 1$, which follows from the binomial theorem. We also have performed most of our calculations via \textit{Mathematica}. 

\section{Proof of Theorem \ref{thm11}}\label{sec2}

We first prove (\ref{eq16}). We start with the following $2$-dissection \cite[(2.1.1)]{humcy}
\begin{align}\label{eq21}
\dfrac{f_2^5}{f_1^2f_4^2} = \dfrac{f_8^5}{f_4^2f_{16}^2}+2q\dfrac{f_{16}^2}{f_8}.
\end{align}
Dividing both sides of (\ref{eq21}) by $f_2f_4$ yields
\begin{align}\label{eq22}
\sum_{n= 0}^\infty B(n)q^n = \dfrac{f_2^4}{f_1^2f_4^3} = \dfrac{f_8^5}{f_2f_4^3f_{16}^2}+2q\dfrac{f_{16}^2}{f_2f_4f_8}.
\end{align}
We consider the terms in (\ref{eq22}) involving $q^{2n+1}$, so that
\begin{align*}
\sum_{n=0}^\infty B(2n+1)q^n = 2\dfrac{f_8^2}{f_1f_2f_4},
\end{align*}
which immediately implies (\ref{eq16}). 

We next prove (\ref{eq17}). We require the $q$-series identities \cite[(1.7.1), (1.5.5)]{hirscb}
\begin{align}
f_1^3 &= \sum_{k=0}^\infty (-1)^k(2k+1)q^{k(k+1)/2},\label{eq23}\\
\dfrac{f_2^2}{f_1} &= \sum_{k=0}^\infty q^{k(k+1)/2}.\label{eq24}
\end{align}
Using (\ref{eq23}) and (\ref{eq24}), we express
\begin{align}
\sum_{n= 0}^\infty B(n)q^n &= \dfrac{f_2^4}{f_1^2f_4^3} \equiv \dfrac{f_{10}}{f_5f_{20}}\cdot \dfrac{f_4^2}{f_2}\cdot f_1^3\nonumber\\
&\equiv \dfrac{f_{10}}{f_5f_{20}}\sum_{k=0}^\infty\sum_{l=0}^\infty (-1)^l(2l+1)q^{k(k+1)+l(l+1)/2}\pmod{5}.\label{eq25}
\end{align}
We consider the equation 
\begin{align*}
k(k+1)+\dfrac{l(l+1)}{2}\equiv 4\pmod{5},
\end{align*}
which is equivalent to
\begin{align}\label{eq26}
2(2k+1)^2+(2l+1)^2\equiv 0\pmod{5}.
\end{align}
Since $-2$ is a quadratic nonresidue modulo $5$, we infer from (\ref{eq26}) that $2k+1\equiv 2l+1\equiv0\pmod{5}$. Thus, by considering the terms of (\ref{eq25}) involving $q^{5n+4}$, we arrive at (\ref{eq17}). \hfill $\square$

\section{Proofs of Theorems \ref{thm13} and \ref{thm14}}\label{sec3}

We obtain in this section the generating function for $B(3n+2)$, which will be needed to deduce Theorems \ref{thm13} and \ref{thm14}.

\begin{lemma}\label{lem31}
We have the identity 
\begin{align*}
	\sum_{n=0}^\infty B(3n+2)q^n = \dfrac{f_2^{12}f_{12}^3}{f_1^6f_4^{10}}.
\end{align*}
\end{lemma}

\begin{proof}
We need the following $3$-dissections \cite[(2.2.1), (2.2.9)]{humcy}
\begin{align}
	\dfrac{f_2^2}{f_1} &= \dfrac{f_6f_9^2}{f_3f_{18}} + q\dfrac{f_{18}^2}{f_9},\label{eq31}\\
	\dfrac{1}{f_1^3} &= \dfrac{f_9^3}{f_3^{10}}\left(\alpha(q^3)^2+3q\dfrac{f_9^3}{f_3}\alpha(q^3)+9q^2\frac{f_9^6}{f_3^2}\right)\label{eq32}
\end{align}
and the identity \cite[Lemma 2.2]{guad1}
\begin{align}\label{eq33}
	\dfrac{f_6^2}{f_3}\alpha(q^4) + 3q\dfrac{f_2f_3^2f_{12}^3}{f_1f_4f_6} = \dfrac{f_2^6}{f_1^3},
\end{align}
We apply (\ref{eq31}) and (\ref{eq32}) on the generating function for $B(n)$ so that
\begin{align}\label{eq34}
	\sum_{n= 0}^\infty B(n)q^n &= \left(\dfrac{f_2^2}{f_1}\right)^2\cdot \dfrac{1}{f_4^3}\nonumber\\
	&= \left(\dfrac{f_6f_9^2}{f_3f_{18}} + q\dfrac{f_{18}^2}{f_9}\right)^2 \cdot \dfrac{f_{36}^3}{f_{12}^{10}}\left(\alpha(q^{12})^2+3q^4\dfrac{f_{36}^3}{f_{12}}\alpha(q^{12})+9q^8\frac{f_{}^6}{f_3^2}\right)
\end{align}
We extract the terms of (\ref{eq34}) involving $q^{3n+2}$. In view of (\ref{eq33}), we see that
\begin{align*}
	\sum_{n= 0}^\infty B(3n+2)q^n &= \dfrac{f_6^4f_{12}^3\alpha(q^4)^2}{f_4^{10}f_3^2}+6q\dfrac{f_2f_3f_6f_{12}^6\alpha(q^4)}{f_1f_4^{11}}+9q^2\dfrac{f_2^2f_3^4f_{12}^9}{f_1^2f_4^{12}f_6^2}\\
	&=\dfrac{f_{12}^3}{f_4^{10}}\left(\dfrac{f_6^2}{f_3}\alpha(q^4) + 3q\dfrac{f_2f_3^2f_{12}^3}{f_1f_4f_6}\right)^2=\dfrac{f_2^{12}f_{12}^3}{f_1^6f_4^{10}}
\end{align*}
as desired.
\end{proof}

\begin{proof}[Proof of Theorem \ref{thm13}]
We first recall Euler's identity \cite[(1.6.1)]{hirscb}
\begin{align}\label{eq35}
	f_1 &= \sum_{k=-\infty}^\infty (-1)^k q^{k(3k+1)/2}
\end{align}
so that from Lemma \ref{lem31}, we have
\begin{align*}
	\dfrac{f_2^{12}f_{12}^3}{f_1^6f_4^9} = \dfrac{f_2^{12}f_{12}^3}{f_1^6f_4^{10}}\cdot f_4&= \left(\sum_{m= 0}^\infty B(3m+2)q^m\right)\left( \sum_{k=-\infty}^\infty (-1)^k q^{2k(3k+1)}\right)\\
	&= \sum_{n= 0}^\infty C(n)q^n,
\end{align*}
where 
\begin{align}\label{eq36}
	C(n) := \sum_{k=-\infty}^\infty (-1)^k B\left(3n+2-6k(3k+1)\right).
\end{align}
Observe that 
\begin{align}\label{eq37}
	\sum_{n= 0}^\infty C(n)q^n = \dfrac{f_2^{12}f_{12}^3}{f_1^6f_4^9}\equiv \dfrac{f_6^4}{f_3^2}\pmod{3},
\end{align}
whose $q$-expansion modulo $3$ contains only terms of the form $q^{3n}$. Thus, by looking at the terms of (\ref{eq37}) involving $q^{3n+1}$ and $q^{3n+2}$ and using (\ref{eq36}), we arrive at the desired congruence.
\end{proof}

\begin{proof}[Proof of Theorem \ref{thm14}]
We extract the terms of (\ref{eq37}) involving $q^{3n}$, so that
\begin{align*}
	\sum_{n= 0}^\infty C(3n)q^n \equiv \dfrac{f_2^4}{f_1^2} \pmod{3}.
\end{align*}
Suppose $p\geq 5$ is a prime with $p\equiv 3\pmod{4}$. Setting $(r,s)=(-2,4)$ in \cite[Theorem 2]{cooph} gives
\begin{align}\label{eq38}
	C\left(3pn+\dfrac{3(p^2-1)}{4}\right)\equiv C\left(\dfrac{3n}{p}\right)\pmod{3},
\end{align}
where we set $C(n) = 0$ if $n$ is not an integer. Replacing $n$ with $pn+r$, where $r\in \{1,\ldots,p-1\}$, in (\ref{eq38}) yields
\begin{align}\label{eq39}
	C\left(3p(pn+r)+\dfrac{3(p^2-1)}{4}\right)\equiv 0\pmod{3}.
\end{align}
Combining (\ref{eq36}) and (\ref{eq39}), we get the desired congruence.
\end{proof}

\section{Proofs of Theorems \ref{thm12}, \ref{thm15}, and \ref{thm16}}\label{sec4}

We first derive in this section the generating function modulo $9$ for $B(3n+1)$, which will be used to prove Theorem \ref{thm15}. We present the following lemma.

\begin{lemma}\label{lem41}
We have the identity
\begin{align*}
	f_{12}^3f_2^3f_3^6 + f_1^2f_4f_6^9 = 2f_1f_2f_3^3f_4^3f_6^2f_{12}^2.
\end{align*}
\end{lemma}

\begin{proof}
We need the identities \cite[(12.18), (12.19), (12.21)]{coop}
\begin{align}
	\dfrac{1}{h}+h &= \dfrac{f_3^3f_4}{qf_1f_{12}^3},\label{eq41}\\
	\dfrac{1}{h}-1+h &= \dfrac{f_4^4f_6^2}{qf_2^2f_{12}^4},\label{eq42}\\
	\dfrac{1}{h}-2+h &= \dfrac{f_1f_4^2f_6^9}{qf_2^3f_3^3f_{12}^6},\label{eq43}
\end{align}
where 
\begin{align*}
	h := h(q) = q\prod_{n=1}^\infty \dfrac{(1-q^{12n-1})(1-q^{12n-11})}{(1-q^{12n-5})(1-q^{12n-7})}
\end{align*}
is the level $12$ analogue of the Rogers–Ramanujan continued fraction. We have that
\begin{align}
	\dfrac{f_3^3f_4}{qf_1f_{12}^3}+\dfrac{f_1f_4^2f_6^9}{qf_2^3f_3^3f_{12}^6} &= \dfrac{1}{h}+h+\dfrac{1}{h}-2+h\nonumber\\
	&= 2\left(\dfrac{1}{h}-1+h\right) = 2\dfrac{f_4^4f_6^2}{qf_2^2f_{12}^4}. \label{eq44}
\end{align}
Multiplying both sides of (\ref{eq44}) by $qf_1f_2^3f_3^3f_{12}^6/f_4$ leads to the desired identity.
\end{proof}

\begin{lemma}\label{lem42}
We have the identity 
\begin{align*}
	\sum_{n=0}^\infty B(3n+1)q^n \equiv 2\dfrac{f_1^2f_2^4}{f_4}\pmod{9}.
\end{align*}
\end{lemma}

\begin{proof}
We look for the terms of (\ref{eq34}) involving $q^{3n+1}$. We obtain 
\begin{align}\label{eq45}	
	\sum_{n= 0}^\infty B(3n+1)q^n &\equiv 2\alpha^2(q^4)\dfrac{f_2f_3f_6f_{12}^3}{f_1f_4^{10}}+3q\alpha(q^4)\dfrac{f_2^2f_3^4f_{12}^6}{f_1^2f_6^2f_4^{11}}\pmod{9}.
\end{align}
Invoking the identities \cite[(22.11.6), (22.6.1)]{hirscb}
\begin{align*}
	\alpha(q) &= \dfrac{f_2^6f_3}{f_1^3f_6^2}+3q \dfrac{f_6^6f_1}{f_3^3f_2^2},\\
	\alpha(q^4) &= \alpha(q) - 6q\dfrac{f_4^2f_{12}^2}{f_2f_6},
\end{align*}
we surmise that
\begin{align}\label{eq46}
	\alpha(q^4) = \dfrac{f_2^6f_3}{f_1^3f_6^2}+3q\left(\dfrac{f_6^6f_1}{f_3^3f_2^2}-2\dfrac{f_4^2f_{12}^2}{f_2f_6}\right)
\end{align}
and
\begin{align}\label{eq47}
	\alpha^2(q^4) \equiv \dfrac{f_2^{12}f_3^2}{f_1^6f_6^4} + 6q\left(\dfrac{f_2^4f_6^4}{f_1^2f_3^2}-2\dfrac{f_2^5f_3f_4^2f_{12}^2}{f_1^3f_6^3}\right)\pmod{9}.
\end{align}
Substituting (\ref{eq46}) and (\ref{eq47}) into (\ref{eq45}), we express
\begin{align}	
	\sum_{n= 0}^\infty B(3n+1)q^n &\equiv 2\dfrac{f_2f_3f_6f_{12}^3}{f_1f_4^{10}}\left[\dfrac{f_2^{12}f_3^2}{f_1^6f_6^4} + 6q\left(\dfrac{f_2^4f_6^4}{f_1^2f_3^2}-2\dfrac{f_2^5f_3f_4^2f_{12}^2}{f_1^3f_6^3}\right)\right]+3q\dfrac{f_2^2f_3^4f_{12}^6}{f_1^2f_6^2f_4^{11}}\cdot \dfrac{f_2^6f_3}{f_1^3f_6^2}\nonumber\\
	&\equiv2\dfrac{f_2^{13}f_3^3f_{12}^3}{f_1^7f_4^{10}f_6^3}+3q\left(\dfrac{f_2^5f_6^5f_{12}^3}{f_1^3f_3f_4^{10}}-2\dfrac{f_2^6f_3^2f_{12}^5}{f_1^4f_4^8f_6^2}+\dfrac{f_2^8f_3^5f_{12}^6}{f_1^5f_4^{11}f_6^4}\right)\pmod{9}.\label{eq48}
\end{align}
Applying Lemma \ref{lem41}, we deduce that
\begin{align}
	\dfrac{f_2^5f_6^5f_{12}^3}{f_1^3f_3f_4^{10}}-2\dfrac{f_2^6f_3^2f_{12}^5}{f_1^4f_4^8f_6^2}&+\dfrac{f_2^8f_3^5f_{12}^6}{f_1^5f_4^{11}f_6^4}\nonumber\\
	&=\dfrac{f_2^5f_{12}^3}{f_1^5f_3f_4^{11}f_6^4}\left(f_1^2f_4f_6^9-2f_1f_2f_3^3f_4^3f_6^2f_{12}^2+f_{12}^3f_2^3f_3^6\right)\nonumber\\
	&=0.\label{eq49}
\end{align}
Hence, combining (\ref{eq48}) and (\ref{eq49}) yields
\begin{align*}
	\sum_{n= 0}^\infty B(3n+1)q^n \equiv 2\dfrac{f_2^{13}f_3^3f_{12}^3}{f_1^7f_4^{10}f_6^3}\equiv 2\dfrac{f_1^2f_2^4}{f_4}\pmod{9}
\end{align*}
as desired.
\end{proof}

\begin{proof}[Proof of Theorem \ref{thm15}]
Using (\ref{eq35}) and Lemma \ref{lem42}, we know that
\begin{align}
	2f_1^2f_2^4 = 2\dfrac{f_1^2f_2^4}{f_4}\cdot f_4&\equiv \left(\sum_{m= 0}^\infty B(3m+1)q^m\right)\left( \sum_{k=-\infty}^\infty (-1)^k q^{2k(3k+1)}\right)\nonumber\\
	&\equiv \sum_{n= 0}^\infty D(n)q^n \pmod{9},\label{eq410}
\end{align}
where 
\begin{align}\label{eq411}
	D(n) := \sum_{k=-\infty}^\infty (-1)^k B\left(3n+1-6k(3k+1)\right).
\end{align}
Let $p\equiv 7, 11\pmod{12}$ be a prime. We use $(r,s) = (2,4)$ in \cite[Theorem 2]{cooph} in (\ref{eq410}) so that
\begin{align}\label{eq412}
	D\left(pn+\dfrac{5(p^2-1)}{12}\right)\equiv p^2D\left(\dfrac{n}{p}\right)\pmod{9},
\end{align}
where $D(n)=0$ if $n$ is not an integer. Replacing $n$ with $pn+r$, where $r\in \{1,\ldots,p-1\}$, in (\ref{eq412}) yields
\begin{align}\label{eq413}
	D\left(p(pn+r)+\dfrac{5(p^2-1)}{12}\right)\equiv 0\pmod{9}.
\end{align}
Combining (\ref{eq411}) and (\ref{eq413}) leads us to the desired congruence.
\end{proof}

We next find the generating function modulo $9$ for $B(9n+7)$, which will be needed to prove Theorems \ref{thm12} and \ref{thm16}.

\begin{lemma}\label{lem43}
We have the identity
\begin{align*}
	\sum_{n=0}^\infty B(9n+7)q^n \equiv -\dfrac{f_1^3f_3f_4^2f_{12}}{f_2f_6}\pmod{9}.
\end{align*}
\end{lemma}

\begin{proof}
We require the following $3$-dissections \cite[(2.2.2), (2.2.3)]{humcy}
\begin{align}
	\dfrac{f_1f_4}{f_2} &= \dfrac{f_3f_{12}f_{18}^5}{f_6^2f_9^2f_{36}^2} - q\dfrac{f_9f_{36}}{f_{18}},\label{eq414}\\
	\dfrac{f_1^2}{f_2} &= \dfrac{f_9^2}{f_{18}} - 2q\dfrac{f_3f_{18}^2}{f_6f_9}.\label{eq415}
\end{align}
Applying (\ref{eq414}) and (\ref{eq415}) on Lemma \ref{lem42} gives
\begin{align}
	\sum_{n=0}^\infty B(3n+1)q^n &\equiv 2\dfrac{f_1^2f_2^4}{f_4}\equiv 2\left(\dfrac{f_1f_4}{f_2}\right)^2\left(\dfrac{f_2^2}{f_4}\right)^3\nonumber\\
	&\equiv 2\left(\dfrac{f_3f_{12}f_{18}^5}{f_6^2f_9^2f_{36}^2} - q\dfrac{f_9f_{36}}{f_{18}}\right)\left(\dfrac{f_{18}^2}{f_{36}} - 2q^2\dfrac{f_6f_{36}^2}{f_{12}f_{18}}\right)^3\pmod{9}.\label{eq416}
\end{align}
We look at the terms of (\ref{eq416}) involving $q^{3n+2}$. We then deduce that
\begin{align}
	\sum_{n=0}^\infty B(9n&+7)q^n \nonumber\\
	&\equiv 2\dfrac{f_3^2f_6^4}{f_{12}}-3\dfrac{f_1^2f_4f_6^{13}}{f_2^3f_3^4f_{12}^4}-3q\dfrac{f_1f_6^4f_{12}^2}{f_3f_4}+2q^2\dfrac{f_2^3f_3^2f_{12}^8}{f_4^3f_6^5}\nonumber\\
	&\equiv q\dfrac{f_1f_6^4f_{12}^2}{f_3f_4}\left(2\dfrac{f_3^3f_4}{qf_1f_{12}^3}-3\dfrac{f_1f_4^2f_6^9}{qf_2^3f_3^3f_{12}^6}-3+2q\dfrac{f_2^3f_3^3f_{12}^6}{f_1f_4^2f_6^9}\right)\pmod{9}.\label{eq417}
\end{align}
We now employ (\ref{eq41}), (\ref{eq42}), (\ref{eq43}), and the identity \cite[(12.20)]{coop}
\begin{align*}
	\dfrac{1}{h}-4+h = \dfrac{f_1^3f_4f_6^2}{qf_2^2f_3f_{12}^3}
\end{align*}
to simplify the right-hand side of (\ref{eq417}). We compute
\begin{align}
	2\dfrac{f_3^3f_4}{qf_1f_{12}^3}&-3\dfrac{f_1f_4^2f_6^9}{qf_2^3f_3^3f_{12}^6}-3+2q\dfrac{f_2^3f_3^3f_{12}^6}{f_1f_4^2f_6^9}\nonumber\\
	&= 2\left(\dfrac{1}{h}+h\right)-3\left(\dfrac{1}{h}-2+h\right)-3+2\left(\dfrac{1}{h}-2+h\right)^{-1}\nonumber\\
	&= -\dfrac{1-4h+h^2}{h}\cdot \dfrac{1-h+h^2}{1-2h+h^2}\nonumber\\
	&= -\left(\dfrac{1}{h}-4+h\right)\left(\dfrac{1}{h}-1+h\right)\left(\dfrac{1}{h}-2+h\right)^{-1}\nonumber\\
	&= -\dfrac{f_1^3f_4f_6^2}{qf_2^2f_3f_{12}^3}\cdot \dfrac{f_4^4f_6^2}{qf_2^2f_{12}^4}\cdot q\dfrac{f_2^3f_3^3f_{12}^6}{f_1f_4^2f_6^9}=-\dfrac{f_1^2f_3^2f_4^3}{qf_2f_6^5f_{12}}.\label{eq418}
\end{align}
We see from (\ref{eq417}) and (\ref{eq418}) that
\begin{align*}
	\sum_{n=0}^\infty B(9n+7)q^n\equiv -q\dfrac{f_1f_6^4f_{12}^2}{f_3f_4}\cdot \dfrac{f_1^2f_3^2f_4^3}{qf_2f_6^5f_{12}} \equiv -\frac{f_1^3f_3f_4^2f_{12}}{f_2f_6}\pmod{9}
\end{align*}
as desired.
\end{proof}

\begin{proof}[Proof of Theorem \ref{thm12}]
We need the following $3$-dissection \cite[(2.2.8)]{humcy}
\begin{align}\label{eq419}
	f_1^3 = \alpha(q^3)f_3 - 3qf_9^3.
\end{align}
We apply (\ref{eq31}) and (\ref{eq419}) on Lemma \ref{lem43} and get
\begin{align}
	\sum_{n=0}^\infty B(9n+7)q^n &\equiv -\frac{f_1^3f_3f_4^2f_{12}}{f_2f_6} \equiv -\dfrac{f_3f_{12}}{f_6}\cdot f_1^3\cdot \dfrac{f_4^2}{f_2}\nonumber\\
	&\equiv -\dfrac{f_3f_{12}}{f_6}(\alpha(q^3)f_3 - 3qf_9^3)\left(\dfrac{f_{12}f_{18}^2}{f_6f_{36}} + q^2\dfrac{f_{36}^2}{f_{18}}\right)\pmod{9}.\label{eq420}
\end{align}
We extract the terms of (\ref{eq420}) so that
\begin{align}\label{eq421}
	\sum_{n=0}^\infty B(27n+16)q^n \equiv \dfrac{f_1f_4}{f_2}\cdot 3f_3^3 \cdot \dfrac{f_4f_6^2}{f_2f_{12}} \equiv 3\dfrac{f_1f_3^3f_4^2f_6^2}{f_2^2f_{12}}\pmod{9},
\end{align}
which immediately implies the desired congruence.
\end{proof}

\begin{proof}[Proof of Theorem \ref{thm16}]
We use the following $q$-series identity \cite[(10.7.7)]{hirscb}
\begin{align*}
	\dfrac{f_2^5}{f_1^2} &= \sum_{k=-\infty}^\infty (-1)^k(3k+1)q^{k(3k+2)}
\end{align*}
on (\ref{eq421}) so that 
\begin{align}
	3\dfrac{f_1f_3^3f_4^7f_6^2}{f_2^4f_{12}}&\equiv 3\dfrac{f_1f_3^3f_4^2f_6^2}{f_2^2f_{12}}\cdot \dfrac{f_4^5}{f_2^2}\nonumber\\
	&\equiv \left(\sum_{m=0}^\infty B(27m+16)q^m\right)\left(\sum_{k=-\infty}^\infty (-1)^k(3k+1)q^{2k(3k+2)}\right)\nonumber\\
	&\equiv \sum_{n= 0}^\infty E(n)q^n \pmod{9}, \label{eq422}
\end{align}
where 
\begin{align}
	E(n) := \sum_{k=-\infty}^\infty (-1)^k(3k+1)B(27n+16-54k(3k+2)).
\end{align}
We now employ (\ref{eq414}) in (\ref{eq422}), yielding 
\begin{align}
	\sum_{n= 0}^\infty E(n)q^n \equiv 3\dfrac{f_1f_3^3f_4^7f_6^2}{f_2^4f_{12}}&\equiv 3f_3^3f_6f_{12}\cdot \dfrac{f_1f_4}{f_2}\nonumber\\
	&\equiv 3f_3^3f_6f_{12}\left(\dfrac{f_3f_{12}f_{18}^5}{f_6^2f_9^2f_{36}^2} - q\dfrac{f_9f_{36}}{f_{18}}\right)\pmod{9}.\label{eq423}
\end{align}
Considering the terms of (\ref{eq423}) involving $q^{3n+2}$, we arrive at $E(3n+2)\equiv 0\pmod{9}$ for $n\geq 0$. Hence, the desired congruence follows from this congruence and (\ref{eq423}).
\end{proof}

\section{Proof of Theorem \ref{thm17}}\label{sec5}

We derive in this section the generating function modulo $7$ for $B(7n+2)$ by using the \textit{Mathematica} package \texttt{RaduRK} created by Smoot \cite{smoot} based from Radu's Ramanujan--Kolberg algorithm \cite{radu}. Before we explain how this algorithm works, we first give a brief background on modular functions on the congruence subgroup 
\begin{align*}
\Gamma_0(N) &:= \left\lbrace\begin{bmatrix}
	a & b\\ c& d
\end{bmatrix} \in \mbox{SL}_2(\mathbb{Z}) : c\equiv 0\pmod N\right\rbrace.
\end{align*}
Recall that a matrix $\gamma:=[\begin{smallmatrix}
a & b\\ c& d
\end{smallmatrix}]\in \Gamma_0(N)$ acts on an element $\tau$ in the extended upper-half plane $\mathbb{H}^\ast := \mathbb{H}\cup \mathbb{Q}\cup \{\infty\}$ via 
\begin{align*}
\gamma\tau = \dfrac{a\tau+b}{c\tau+d}.
\end{align*}
We define the cusps of $\Gamma_0(N)$ as the equivalence classes of $\mathbb{Q}\cup \{\infty\}$ under this action. We define a modular function on $\Gamma_0(N)$ as a meromorphic function $f:\mathbb{H}\rightarrow\mathbb{C}$ such that $f(\gamma\tau)=f(\tau)$ for all $\gamma\in\Gamma_0(N)$ and for every cusp $a/c$ of $\Gamma_0(N)$ and $\gamma\in\mbox{SL}_2(\mathbb{Z})$ with $\gamma(\infty)=a/c$, we have the $q$-expansion given by
\begin{align*}
f(\gamma\tau) = \sum_{n= n_0}^\infty a_nq^{n\gcd(c^2,N)/N}
\end{align*}
for some integer $n_0$ with $a_{n_0}\neq 0$, where $q:=e^{2\pi i\tau}$. The integer $n_0$ is called the order of $f$ at $a/c$, and we call $a/c$ a zero (respectively, a pole) of $f(\tau)$ if its order at $a/c$ is positive (respectively, negative).

We define an $\eta$-quotient as a product of the form 
\begin{align*}
f(\tau) = \prod_{\delta\mid N} \eta^{r_{\delta}}(\delta\tau)
\end{align*}
for some indexed set $\{r_\delta\in\mathbb{Z} : \delta\mid N\}$, where $\eta(\tau) := q^{1/24}f_1$ is the Dedekind eta function. One may impose conditions on $r_\delta$ that will make a given $\eta$-quotient modular on $\Gamma_0(N)$ and then compute its orders at the cusps of $\Gamma_0(N)$; we refer the reader to \cite[Theorem 1.64]{ono} and \cite[Theorem 1.65]{ono} for the precise statements. 

Let $\mathcal{M}^\infty(\Gamma_0(N))$ be the algebra of all modular functions on $\Gamma_0(N)$ with a pole only at $\infty$ and $\mathcal{E}^\infty(N)$ be its subalgebra comprising all $\eta$-quotients on $\Gamma_0(N)$. For $N\geq 2$, Radu \cite[Section 2]{radu} explicitly described a basis for the algebra $\langle \mathcal{E}^\infty(N)\rangle_{\mathbb{Q}}$ as a finitely generated $\mathbb{Q}[t]$-module for some $t\in \mathcal{M}^\infty(\Gamma_0(N))$. 

Given positive integers $M, N, m$, and $j$, where $N\geq 2$ and $0\leq j < m$, and a sequence $r=(r_\delta)_{\delta\mid M}$ of integers indexed by the positive divisors of $M$, Radu's algorithm takes the generating function
\begin{align*}
\sum_{n=0}^\infty a(n)q^n = \prod_{\delta\mid N} f_{\delta}^{r_\delta}
\end{align*}
as an input and checks if there exist an $\alpha\in\mathbb{Q}$, an $\eta$-quotient $f(\tau)$ on $\Gamma_0(N)$, and a set $P_{m,r}(j)\subseteq \{0,1,\ldots,m-1\}$ uniquely determined by $m, r$, and $j$ such that
\begin{align*}
q^\alpha f(\tau)\prod_{j'\in P_{m,r}(j)}\sum_{n= 0}^\infty a(mn+j')q^n
\end{align*}
is a modular function on $\Gamma_0(N)$. Finding the minimum value of $N$ satisfying this property is related to the $\Delta^\ast$ criterion (see \cite{radu} for more details) and can be obtained by calling the command {\tt minN[M,r,m,j]}. When such an $N$ (or a multiple of it) is found, we can now write
\begin{align}\label{eq51}
q^\alpha f(\tau)\prod_{j'\in P_{m,r}(j)}\sum_{n= 0}^\infty a(mn+j')q^n = \sum_g g p_g(t)
\end{align}
where $g$ runs all over the elements of an algebra basis for $\langle \mathcal{E}^\infty(N)\rangle_{\mathbb{Q}}$ viewed as a finitely generated $\mathbb{Q}[t]$-module for some $t\in \mathcal{M}^\infty(\Gamma_0(N))$ and $p_g(X)$ are polynomials in $X$ with integer coefficients. The identity (\ref{eq51}) then yields the product of the generating functions for $a(mn+j')$ when $j'$ runs all over the elements of $P_{m,r}(j)$.

We now use Radu's algorithm to deduce the generating function modulo $7$ for $B(7n+2)$, as this will be needed to prove Theorem \ref{thm17}.

\begin{lemma}\label{lem51}
We have the identity 
\begin{align*}
	\sum_{n=0}^\infty B(7n+2)q^n \equiv \dfrac{f_1f_2^2f_4^2f_{14}^2}{f_7f_{28}}\pmod{7}.
\end{align*}
\end{lemma}

\begin{proof}
Looking at the generating function for $B(n)$, we set $(M,m,j) = (4,7,2)$ and $r=(-2,4,-3)$. We run the command {\tt minN[4,\{-2,4,-3\},7,2]}, which outputs $N=28$. This means that we work on the congruence subgroup $\Gamma_0(28)$ to derive an identity of the form (\ref{eq51}).

We first construct an algebra basis for $\langle \mathcal{E}^\infty(28)\rangle_\infty$. Since the corresponding modular curve $X_0(28) := \Gamma_0(28)\setminus \mathbb{H}^\ast$ has genus $2$, the Weierstrass gap theorem \cite[Theorem 1.1]{paulerd2} dictates that any element of $\mathcal{M}^\infty(\Gamma_0(28))$ must have pole order of at most $-3$. In view of a refinement of Newman's conjecture due to Paule and Radu \cite[Conjecture 9.4]{paulerd1}, a sufficient algebra basis for $\langle \mathcal{E}^\infty(28)\rangle_\infty$ must contain $\eta$-quotients whose orders at $\infty$ are $-3, -4$, and $-5$. To find such $\eta$-quotients, we just run {\tt e28 := etaGenerators[28];} and set 
\begin{align*}
	X &:= {\tt e28[[1]]} = \dfrac{\eta^4(4\tau)\eta^2(14\tau)}{\eta^2(2\tau)\eta^4(28\tau)} = \dfrac{f_4^4f_{14}^2}{q^3f_2^2f_{28}^4},\\
	Y &:= {\tt e28[[2]]} = \dfrac{\eta(2\tau)\eta^2(4\tau)\eta^5(14\tau)}{\eta(\tau)\eta(7\tau)\eta^6(28\tau)} = \dfrac{f_2f_4^2f_{14}^5}{q^4f_1f_7f_{28}^6},\\
	Z &:= {\tt e28[[6]]} = \dfrac{\eta(2\tau)\eta(4\tau)\eta^5(14\tau)}{\eta^7(28\tau)} = \dfrac{f_2f_4f_{14}^5}{q^5f_{28}^7}.
\end{align*}
We now define the aforementioned basis as 
\begin{align*}
	{\tt \{e28[[1]],\{1,e28[[2]],e28[[6]]\}\}}
\end{align*}  
and then run the command 
\begin{align}\label{eq52}
	{\tt RKMan[28,4,\{-2,4,-3\},7,2,\{e28[[1]],\{1,e28[[2]],e28[[6]]\}\}]}.
\end{align}
We then obtain the following output, as shown in Table \ref{tab51} below. 

\begin{table}[ht]
	\caption{Output of the command (\ref{eq52})}\label{tab51}%
	\begin{center}
		\begin{tabular}{@{}l|l@{}}
			\toprule
			$P_{m,r}(j):$ & $\{2\}$\\
			\midrule
			$q^{\alpha}f(\tau):$ & $\dfrac{f_1^6f_4^{19}f_{14}^{13}}{q^{20}f_2^{11}f_{28}^{26}}$\\
			\midrule
			$t:$ & $\dfrac{f_4^4f_{14}^2}{q^3f_2^2f_{28}^4}$\\
			\midrule
			$\mbox{AB}$: & $\left\{1,\dfrac{f_2f_4^2f_{14}^5}{q^4f_1f_7f_{28}^6},\dfrac{f_2f_4f_{14}^5}{q^5f_{28}^7}\right\}$\\
			\midrule
			$\{p_g(t) : g\in\mbox{AB}\}$: & $\begin{aligned}
				&\{-2401t-5145t^2+6860t^3+882t^4-175t^5-21t^6,\\
				& 2401-7154t^2-294t^3+189t^4+14t^5,\\
				& 3430t-735t^3-42t^4+t^5\}
			\end{aligned}$\\
			\midrule
			Common Factor : & none\\
			\bottomrule
		\end{tabular}
	\end{center}
\end{table}
In view of (\ref{eq51}), we infer from Table \ref{tab51} that
\begin{align}\label{eq53}
	\dfrac{f_1^6f_4^{19}f_{14}^{13}}{q^{20}f_2^{11}f_{28}^{26}}&\sum_{n=0}^\infty B(7n+2)q^n = (-2401X-5145X^2+6860X^3+882X^4\nonumber\\
	&-175X^5-21X^6)+Y(2401-7154X^2-294X^3+189X^4+14X^5)\nonumber\\
	&+Z(3430X-735X^3-42X^4+X^5).
\end{align}
Taking both sides of (\ref{eq53}) modulo $7$ yields
\begin{align*}
	\sum_{n=0}^\infty B(7n+2)q^n \equiv \dfrac{q^{20}f_2^{11}f_{28}^{26}}{f_1^6f_4^{19}f_{14}^{13}}ZX^5 \equiv \dfrac{f_1f_2^2f_4^2f_{14}^2}{f_7f_{28}}\pmod{7}
\end{align*}
as desired.
\end{proof}

\begin{proof}[Proof of Theorem \ref{thm17}]
We first show (\ref{eq18}). We begin with the following $q$-series identity \cite[(10.7.3)]{hirscb}
\begin{align*}
	\dfrac{f_1^5}{f_2^2} = \sum_{k=-\infty}^\infty (6k+1)q^{k(3k+1)/2},
\end{align*}
so that from Lemma \ref{lem51}, we have
\begin{align}
	\dfrac{f_1f_{14}^3}{f_7f_{28}}\equiv\dfrac{f_1f_2^2f_4^2f_{14}^2}{f_7f_{28}}\cdot \dfrac{f_2^5}{f_4^2} &\equiv \left(\sum_{m= 0}^\infty B(7m+2)q^m\right)\left( \sum_{k=-\infty}^\infty (6k+1)q^{k(3k+1)}\right)\nonumber\\
	&\equiv \sum_{n= 0}^\infty F(n)q^n\pmod{7},\label{eq54}
\end{align}
where 
\begin{align}\label{eq55}
	F(n) := \sum_{k=-\infty}^\infty (6k+1)B\left(7n+2-7k(3k+1)\right).
\end{align}
We next employ the $7$-dissection of $f_1$ \cite[(10.5.1)]{hirscb} given by 
\begin{align*}
	f_1 = f_{49}\left(A_0(q^7) - qA_1(q^7) -q^2+q^5A_5(q^7)\right)
\end{align*}
for some $A_0(q), A_1(q), A_5(q)\in\mathbb{Z}[[q]]$ in (\ref{eq54}) so that
\begin{align}\label{eq56}
	\sum_{n= 0}^\infty F(n)q^n\equiv \dfrac{f_{14}^3f_{49}}{f_7f_{28}}\left(A_0(q^7) - qA_1(q^7) -q^2+q^5A_5(q^7)\right)\pmod{7}.
\end{align}
We read the terms of (\ref{eq56}) involving $q^{7n+j}$, where $j\in \{3,4,6\}$. We obtain $F(7n+j)\equiv 0\pmod{7}$ for all $n\geq 0$, and combining this with (\ref{eq55}), we arrive at (\ref{eq18}). 

We next show (\ref{eq19}). We consider the terms of (\ref{eq56}) involving $q^{7n+2}$, so that
\begin{align}\label{eq57}
	\sum_{n= 0}^\infty F(7n+2)q^n\equiv \dfrac{f_2^3f_7}{f_1f_4}\pmod{7}. 
\end{align}
We replace $q$ with $-q$ in (\ref{eq35}), obtaining
\begin{align*}
	\dfrac{f_2^3}{f_1f_4} = \sum_{k=-\infty}^\infty (-1)^{k(k+1)/2} q^{k(3k+1)/2}.
\end{align*}
Plugging this identity into (\ref{eq57}) yields
\begin{align}\label{eq58}
	\sum_{n= 0}^\infty F(7n+2)q^n\equiv f_7\sum_{k=-\infty}^\infty (-1)^{k(k+1)/2} q^{k(3k+1)/2}\pmod{7}.
\end{align}
Note that $k(3k+1)/2\equiv 0,1,2,5\pmod{7}$. Thus, looking at the terms of (\ref{eq58}) involving $q^{7n+j}$, where $j\in \{3,4,6\}$, we deduce that $F(7(7n+j)+2)\equiv 0\pmod{7}$ for all $n\geq 0$. Hence, (\ref{eq19}) follows from this congruence and (\ref{eq55}).
\end{proof}

\end{document}